\documentclass[a4paper,11pt,reqno]{amsart}

\def\ifundefined#1{\expandafter\ifx\csname#1\endcsname\relax} \ifundefined{preambleloaded}
\usepackage[utf8]{inputenc}
\usepackage{geometry}

\usepackage{amsmath}
\usepackage{amsthm}
\usepackage{amssymb}
\usepackage{amsfonts}
\usepackage{mathtools} % for example for mathrlap
\usepackage{mathrsfs} % for \mathscr
\usepackage{mleftright} % remove extra spacing before \left and after \right
\mleftright % define \left to \mleft and \right to \mright (tighter spacing)

\usepackage[x11names]{xcolor}

\usepackage{tikz-cd}

\usepackage[breaklinks,colorlinks=true,linkcolor=blue,citecolor=green!60!black,backref]{hyperref} % clickable references
\hypersetup{pdftex}
\usepackage{hypcap}
\usepackage[capitalize,nameinlink]{cleveref} % easier references
\crefname{equation}{}{}
\usepackage{autonum}

\usepackage[shortlabels]{enumitem}

\usepackage{csquotes}
\MakeOuterQuote{"}

\makeatletter
\def\@cite#1#2{[{#1\if@tempswa ,\penalty400~#2\fi}]}% NEW
\makeatother

\usetikzlibrary{arrows}
\usetikzlibrary{babel}
\tikzcdset{every label/.append style = {font = \small}}
\tikzset{
    arrowdecor/.style={anchor=south, rotate=90, inner sep=.5mm}
}
\usepackage{rotating}

\newtheorem{thm}{\iflanguage{ngerman}{Satz}{Theorem}}[section]
\newtheorem*{thm*}{\iflanguage{ngerman}{Satz}{Theorem}}
\Crefname{thm}{\iflanguage{ngerman}{Satz}{Theorem}}{\iflanguage{ngerman}{Sätze}{Theorems}}

\newtheorem*{conj*}{\iflanguage{ngerman}{Vermutung}{Conjecture}}
\Crefname{conj}{\iflanguage{ngerman}{Vermutung}{Conjecture}}{\iflanguage{ngerman}{Vermutungen}{Conjectures}}

\newtheorem*{lem*}{Lemma}
\Crefname{lem}{Lemma}{Lemmas}

\newtheorem{prop}[thm]{Proposition}
\newtheorem*{prop*}{Proposition}
\Crefname{prop}{Proposition}{\iflanguage{ngerman}{Propositionen}{Propositions}}

\newtheorem{cor}[thm]{\iflanguage{ngerman}{Korollar}{Corollary}}
\Crefname{cor}{\iflanguage{ngerman}{Korollar}{Corollary}}{\iflanguage{ngerman}{Korollare}{Corollaries}}

\theoremstyle{definition}
\newtheorem{defn}[thm]{Definition}
\Crefname{defn}{Definition}{\iflanguage{ngerman}{Definitionen}{Definitions}}

\newtheorem{defn-prop}[thm]{Definition/Proposition}

\theoremstyle{remark}

\Crefname{rmk}{\iflanguage{ngerman}{Anmerkung}{Remark}}{\iflanguage{ngerman}{Anmerkungen}{Remarks}}

\Crefname{ex}{\iflanguage{ngerman}{Beispiel}{Example}}{\iflanguage{ngerman}{Beispiele}{Examples}}

\newtheoremstyle{claim}
  {0.5\topsep} % Space above
  {0\topsep} % Space below
  {} % Body font
  {} % Indent amount
  {\itshape} % Theorem head font
  {:} % Punctuation after theorem head
  {.5em} % Space after theorem head
  {} % Theorem head spec (can be left empty, meaning `normal')
\theoremstyle{claim}

\newtheorem*{claim*}{\iflanguage{ngerman}{Behauptung}{Claim}}
\Crefname{claim}{\iflanguage{ngerman}{Behauptung}{Claim}}{\iflanguage{ngerman}{Behauptungen}{Claims}}
\counterwithin*{claim}{thm}



\newtheoremstyle{Normal}{}{}{}{}{\bfseries}{:}{.5em}{}
\theoremstyle{Normal}

\counterwithin*{problem}{section}
\makeatletter
\let\oldendproof\endproof
\def\endproof{\oldendproof\aftergroup\@afterindentfalse\aftergroup\@afterheading}
\makeatother
\Crefname{cd}{\iflanguage{ngerman}{Diagramm}{Diagram}}{\iflanguage{ngerman}{Diagramme}{Diagrams}}
\creflabelformat{cd}{(#2\thesection.#1#3)}
\usepackage{aliascnt}
\newaliascnt{cd}{equation}

\makeatletter

\newenvironment{cd*}{%
	$$
	\begin{tikzcd}[sep=large]
}{%
	\end{tikzcd}
	$$\@ignoretrue
}
\makeatother
\makeatletter
\g@addto@macro\bfseries{\boldmath}
\makeatother

\newcommand{\noop}[1]{}

\Crefformat{item}{(#2#1#3)}

\AtBeginDocument{
    \renewcommand{\phi}{\varphi}
    \renewcommand{\epsilon}{\varepsilon}
}

\DeclareMathSymbol{\colonrel}{\mathrel}{operators}{"3A}
\DeclareMathSymbol{:}{\mathpunct}{operators}{"3A}
\expandafter\let\expandafter\originall\csname\encodingdefault\string\l\endcsname
\DeclareRobustCommand*\l
	{\ifmmode\ell\else\expandafter\originall\fi}

\newcommand{\acts}{\raisebox{-0.155em}{\scalebox{0.88}{ \rotatebox[origin=lb]{90}{$\curvearrowleft$}}}\,}

\usepackage{accents}
\DeclareFontFamily{U}{matha}{\hyphenchar\font45}
\DeclareFontShape{U}{matha}{m}{n}{
    <5> <6> <7> <8> <9> <10> gen * matha
    <10.95> matha10 <12> <14.4> <17.28> <20.74> <24.88> matha12
}{}
\DeclareSymbolFont{matha}{U}{matha}{m}{n}
\DeclareFontSubstitution{U}{matha}{m}{n}
\DeclareMathSymbol{\tinsubset}{3}{matha}{"80}
\let\oldbigwedge\bigwedge
\renewcommand{\bigwedge}{\textstyle\oldbigwedge}

\newcommand{\Ric}{\mathrm{Ric}}

\DeclareMathOperator{\Spec}{Spec}

\DeclareMathOperator{\chr}{char}

\DeclareMathOperator{\rk}{rk}

\expandafter\let\expandafter\originald\csname\encodingdefault\string\d\endcsname
\DeclareRobustCommand*\d
	{\ifmmode\mathop{}\!d\else\expandafter\originald\fi}

\newcommand{\p}{\mathrlap{~~\text{.}}}
\expandafter\let\expandafter\originalc\csname\encodingdefault\string\c\endcsname
\DeclareRobustCommand*\c
	{\ifmmode\mathrlap{~~\text{,}}\else\expandafter\originalc\fi}

\makeatletter
\newcommand{\ovl}[1]{%
  \overline{\raisebox{0pt}[\dimexpr\height+0.6pt\relax]{\m@th$#1$}}%
}
\makeatother

\makeatletter
\let\save@mathaccent\mathaccent
\newcommand*\if@single[3]{%
    \setbox0\hbox{${\mathaccent"0362{#1}}^H$}%
    \setbox2\hbox{${\mathaccent"0362{\kern0pt#1}}^H$}%
    \ifdim\ht0=\ht2 #3\else #2\fi
}
\newcommand*\rel@kern[1]{\kern#1\dimexpr\macc@kerna}
\newcommand*\widebar[1]{\@ifnextchar^{{\wide@bar{#1}{0}}}{\wide@bar{#1}{1}}}
\newcommand*\wide@bar[2]{\if@single{#1}{\wide@bar@{#1}{#2}{1}}{\wide@bar@{#1}{#2}{2}}}
\newcommand*\wide@bar@[3]{%
    \begingroup
    \def\mathaccent##1##2{%
        \let\mathaccent\save@mathaccent
        \if#32 \let\macc@nucleus\first@char \fi
        \setbox\z@\hbox{$\macc@style{\macc@nucleus}_{}$}%
        \setbox\tw@\hbox{$\macc@style{\macc@nucleus}{}_{}$}%
        \dimen@\wd\tw@
        \advance\dimen@-\wd\z@
        \divide\dimen@ 3
        \@tempdima\wd\tw@
        \advance\@tempdima-\scriptspace
        \divide\@tempdima 10
        \advance\dimen@-\@tempdima
        \ifdim\dimen@>\z@ \dimen@0pt\fi
        \rel@kern{0.6}\kern-\dimen@
        \if#31
        \overline{\rel@kern{-0.6}\kern\dimen@\macc@nucleus\rel@kern{0.4}\kern\dimen@}%
        \advance\dimen@0.4\dimexpr\macc@kerna
        \let\final@kern#2%
        \ifdim\dimen@<\z@ \let\final@kern1\fi
        \if\final@kern1 \kern-\dimen@\fi
        \else
        \overline{\rel@kern{-0.6}\kern\dimen@#1}%
        \fi
    }%
    \macc@depth\@ne
    \let\math@bgroup\@empty \let\math@egroup\macc@set@skewchar
    \mathsurround\z@ \frozen@everymath{\mathgroup\macc@group\relax}%
    \macc@set@skewchar\relax
    \let\mathaccentV\macc@nested@a
    \if#31
    \macc@nested@a\relax111{#1}%
    \else
    \def\gobble@till@marker##1\endmarker{}%
    \futurelet\first@char\gobble@till@marker#1\endmarker
    \ifcat\noexpand\first@char A\else
    \def\first@char{}%
    \fi
    \macc@nested@a\relax111{\first@char}%
    \fi
    \endgroup
}
\makeatother
\def\semicolon{;}
\def\applytolist#1{
    \expandafter\def\csname multi#1\endcsname##1{
        \def\multiack{##1}\ifx\multiack\semicolon
            \def\next{\relax}
        \else
            \csname #1\endcsname{##1}
            \def\next{\csname multi#1\endcsname}
        \fi
        \next}
    \csname multi#1\endcsname}

\def\calchar#1{\expandafter\newcommand\csname c#1\endcsname{{\mathcal #1}}}
\applytolist{calchar}QWERTYUIOPLKJHGFDSAZXCVBNM;
\def\bbchar#1{\expandafter\newcommand\csname bb#1\endcsname{{\mathbb #1}}}
\applytolist{bbchar}QWERTYUIOPLKJHGFDSAZXCVBNM;
\def\bfchar#1{\expandafter\newcommand\csname bf#1\endcsname{{\mathbf #1}}}
\applytolist{bfchar}qwertyuiopasdfghjklzxcvbnmQWERTYUIOPLKJHGFDSAZXCVBNM;
\def\sfc#1{\expandafter\newcommand\csname s#1\endcsname{{\sf #1}}}
\applytolist{sfc}QWERTYUIOPLKJHGFDSAZXCVBNM;
\def\fchar#1{\expandafter\newcommand\csname f#1\endcsname{{\mathfrak #1}}}
\applytolist{fchar}qwertyuopasdfghjkzxcvbnmQWERTYUIOPLKJHGFDSAZXCVBNM;
\def\scchar#1{\expandafter\newcommand\csname sc#1\endcsname{{\mathscr #1}}}
\applytolist{scchar}QWERTYUIOPLKJHGFDSAZXCVBNM;

\else
\fi

\DeclareMathOperator{\Lie}{Lie}

\DeclareMathOperator{\Fix}{Fix}

\DeclareMathOperator{\ch}{ch}
\DeclareMathOperator{\td}{td}

\DeclareMathOperator{\xxMA}{MA}

\DeclareMathSymbol{\colonrel}{\mathrel}{operators}{"3A}
\DeclareMathSymbol{:}{\mathpunct}{operators}{"3A}

\DeclareRobustCommand*{\MA}[1]{\ifmmode\xxMA#1\else{Monge--Ampère}\fi}
\def\DH/{Duistermaat--Heckman}

\usepackage{bbm}
\renewcommand{\bar}{\widebar}
\renewcommand{\tilde}{\widetilde}

\newcommand{\Oalg}{\mathcal{O}^{\mathrm{alg}}}

\usepackage{microtype}
\counterwithin{equation}{section}

\usepackage{tgpagella}
\newcounter{autoItemCounter}

\newlist{xxthmlist}{enumerate}{1}
\setlist[xxthmlist]{
    label=(\roman{xxthmlisti})\protect\label{item:\theautoItemCounter\thexxthmlisti},
    ref=(\roman{xxthmlisti})
}

\usepackage{thmtools}
\usepackage{thm-restate}

\title{Shrinking gradient Kähler--Ricci solitons are simply-connected}
\author{Carlos Esparza}
\address{Department of Mathematics, University of California Berkeley, CA 94720, USA}
\email{esparzac@berkeley.edu}
\date{\today}

\begin{document}
\begin{abstract}
    We show that smooth polarized Fano fibrations have no nontrivial finite covers.
    Using results by Sun--Zhang and Wylie, it follows that shrinking Kähler--Ricci solitons are simply-connected.
\end{abstract}
\maketitle

\section{Introduction}
A shrinking Kähler--Ricci soliton is a complex manifold $M$ with a complete Kähler metric $\omega$ and $f \in C^\infty(M)$ such that
\[
    \Ric(\omega) + i \partial\bar{\partial} f = \omega
    \p
\]
Sun--Zhang have shown \cite[Thm.~3.1]{SunZhang} that a manifold $M$ admitting such a metric is a quasi-projective variety with additional structure, called a \emph{polarized Fano fibration}:

\begin{defn}[{\cite[Def.~2.5]{SunZhang}}]
    A smooth complex variety $M$ is a Fano fibration if its algebra of regular functions $R \coloneq \Oalg_M(M)$ is finitely generated, the canonical map $\pi: M \to \Spec R$ is a proper surjection and $K_M^{-1}$ is $\pi$-ample.%
    \footnote{
        The definition in \cite{SunZhang} is worded slightly differently. It is easy to see that the affine variety $Y$ in their definition is necessarily $\Spec R$. Furthermore, the ring of regular functions of a smooth (and thus normal) variety is automatically normal.
    }

    A smooth \emph{polarized} Fano fibration consists of a Fano fibration $M$ equipped with the algebraic action of a torus $T$ and an element $X \in \Lie T$ such that under the weight decomposition
    \[
        R = \bigoplus_{\alpha \in (\Lie T)^*} R_\alpha
    \]
    of $T \acts R$ we have $\langle \alpha, X \rangle > 0$ for all non-zero $\alpha$ with $R_\alpha \neq 0$, and $T \acts \Spec R$ only has one fixed point (this last condition is equivalent to $R_0 = \bbC$).
\end{defn}

The purpose of this note is to establish
\begin{prop} \label{thm:pff}
    Let $(M, T, X)$ be a smooth polarized Fano fibration and let $p: \tilde{M} \to M$ be a finite (étale) cover. Then $\tilde{M} = M$.
\end{prop}
It is known that (Riemannian) shrinking gradient Ricci solitons have finite fundamental group \cite{Wylie}, so we conclude:
\begin{cor} \label{thm:main}
    Any shrinking gradient Kähler--Ricci soliton is simply-connected.
\end{cor}

Our proof is inspired by the classical proof that Fano manifolds have no nontrivial finite covers:
If $\tilde{M} \to M$ is a degree $m$ cover of a Fano manifold then $\tilde{M}$ is also Fano.
By the Hirzebruch--Riemann--Roch formula, the holomorphic Euler characteristic is multiplicative along covers, meaning $\chi(\tilde{M}, \cO) = m \chi(M, \cO)$.
But the Hodge numbers $h^{0, p}(M), p \geq 1$ vanish for any Fano manifold by Serre duality and Kodaira vanishing so $\chi(M, \cO) = \chi(\tilde{M}, \cO) = 1$. One concludes that $m = 1$.

To adapt this proof to the noncompact setting of polarized Fano fibrations, we use Wu's equivariant HRR formula \cite[Thm.~3.14]{Wu}.
Sun--Zhang have independently found a different proof of \cref{thm:main}, which instead uses Kähler reduction and results from birational geometry. It will appear in a future version of \cite{SunZhang}.

\subsection*{Acknowledgements.}
I would like to thank Song Sun for his encouragement and guidance, and him and Junsheng Zhang for communicating their proof to me.
I am grateful to thank Charles Cifarelli for conversations that led me to this question and to Zachary Stier for helpful discussions.

\section{Proof of \cref{thm:pff}}
Let $m$ be the order of the cover. We can replace $T$ be its $m$-fold cover in every direction (i.e.\ we replace $\bbR^k / \bbZ^k$ by $\bbR^k / m \bbZ^k$), so that $T$ acts on $\tilde{M}$ and $p$ is $T$-equivariant.

Since $K_M^{-1}$ is relatively ample over a map $M \to Y$ to an affine variety, there exist integers $m, N_1, N_2 > 0$, linear diagonal $T$-actions on $\bbC^{N_1}$ and $\bbP^{N_2}$, and a $T$-equivariant embedding $(j_1, j_2): M \to \bbC^{N_1} \times \bbP^{N_2}$ such that $j_2^*\cO(1) = m K_M^{-1}$ (cf.\ \cite[Def.~2.5]{SunZhang}).
From this embedding we obtain a hermitian metric $h \coloneq e^{-|j_1|^2/2} j_2^* h_\mathrm{FS}^{1/m}$ on $K_M^{-1}$ with positive curvature form $\omega \coloneq (\omega_\mathrm{Euc} \times \omega_\mathrm{FS})|_M$.
By the definition of the Reeb vector field $X$, we have $w_i(X) > 0$ for all weights $w_i \in (\Lie T)^*$ of the action on $\bbC^{N_1}$.
The actions of $T$ on $\bbC^{N_1}$ and $\bbP^{N_2}$ are both Hamiltonian, and the $\omega_\mathrm{Euc}$-moment map for the first action is given by $\sum_i w_i |z_i|^2$. Therefore the action $T \acts M$ is Hamiltonian for $\omega$, and its moment map $\mu: M \to \Lie T$ has the property that $\langle \mu, X \rangle$ is proper (for any Reeb vector field $X$).
Since $p$ is finite, this properness is inherited by $\mu \circ p : \tilde{M} \to \Lie T$, which is a moment map for $T \acts \tilde{M}$.
Thus, both actions $T \acts M$ and $T \acts \tilde{M}$ satisfy Assumption~2.15 in \cite{Wu}, if we pick an action chamber $C$ that contains a Reeb vector field. 

Let $Z = \Fix_T(M)$, with connected components $\{Z_\alpha\}_{\alpha \in F}$.
The covering map $p$ will restrict to (not necessarily connected) covering maps $\tilde{Z}_\alpha \coloneq p^{-1}(Z_\alpha) \to Z_\alpha$, and $\Fix_T(\tilde{M}) = \bigcup_{\alpha \in F} \tilde{Z}_\alpha$.
Applying \cite[Thm.~3.14]{Wu} to the trivial line bundle $\cO$ on $M$ we obtain
\begin{align} \label{eq:local}
    \sum_q (-1)^q \chr H^q(M, \cO) &= \sum_{\alpha \in F} \int_{Z_\alpha} \Omega_{Z_\alpha} \\
    \text{where}\qquad
    \Omega_{Z_\alpha} &= (-1)^{n - n_\alpha - \rk \cN_{Z_\alpha / M}^C} \ch_T \left( \frac{\det(\cN_{Z_\alpha/M}^{-C})}{\det\left(1 - (\cN_{Z_\alpha/M}^C)^*\right) \det\left(1 - \cN_{Z_\alpha/M}^{-C}\right)} \right) \td(Z_\alpha)
\end{align}
and an analogous formula for $\tilde{M}$ and $\{\tilde{Z}_\alpha\}$.
Here $\cN_{Z_\alpha / M}^{\pm C} \leq \cN_{Z_\alpha / M}$ is the subbundle on which elements of $C \subseteq \Lie T$ act with positive/negative weight and $n_\alpha = \dim Z_\alpha = \dim \tilde{Z}_\alpha$.
Evidently $\cN_{\tilde{Z}_\alpha / \tilde{M}} \cong p^* \cN_{Z_\alpha / M}^C$, $T$-equivariantly, so the same is true for the subbundles $\cN^{\pm C}_{Z_\alpha / M}$.

The specific form of the equivariant cohomology classes $\Omega_{Z_\alpha}$ is unimportant.
The crucial point is that the Todd class and the equivariant Chern character $\ch_T$ are functorial along $T$-equivariant maps, so $\Omega_{\tilde{Z}_\alpha} = p^* \Omega_{Z_\alpha}$.
Thus the RHS of \cref{eq:local} is \emph{multiplicative} along the cover $\tilde{M} \to M$, meaning
\[
    \sum_{\alpha \in F} \int_{\tilde{Z}_\alpha} \Omega_{\tilde{Z}_\alpha}
        = m \sum_{\alpha \in F} \int_{Z_\alpha} \Omega_{Z_\alpha}
    \p
\]
Consequently the same is true of the LHS, in particular the constant term of the character on the LHS, i.e.\ 
\[
    \chi^T(M) \coloneq \sum_q (-1)^q \dim H^q(M, \cO)^T
        = \sum_q (-1)^q \dim H^{n, q}(M, K_M^{-1})^T
\]
is multiplicative.

However, since $h$ is a smooth positive metric on $K_M^{-1}$, the cohomology groups $H^{n,q}(N, K_N^{-1})$ vanish for $q \geq 1$ and $N = M, \tilde{M}$.
This is a direct application of \cite[Thm.~VIII.5.6]{agbook}, since $|j_2|^2$ and $|j_2 \circ p|^2$ are weakly plurisubharmonic exhaustion functions on $M$ and $\tilde{M}$ respectively.
Furthermore $H^0(N, \cO)^T = H^0(N, \cO)^{T_\bbC} = \bbC$ because any $J X$-invariant holomorphic function $f$ on $N = M, \tilde{M}$ is constant. This is because $f$ is invariant under the gradient flow of the Hamiltonian $H = \langle \mu, X \rangle$, and since $H$ is Morse--Bott with even indices, a generic gradient flow line approaches the (compact and connected) minimizing submanifold of $H$, on which $f$ must be constant.
Thus $\chi^T(M) = \chi^T(\tilde{M}) = 1$, contradicting the multiplicativity of $\chi^T$ unless $m = 1$.
\hfill $\square$

\bibliographystyle{amsalpha}
\bibliography{complex_bib.bib}

\providecommand{\bysame}{\leavevmode\hbox to3em{\hrulefill}\thinspace}
\providecommand{\MR}{\relax\ifhmode\unskip\space\fi MR }
% \MRhref is called by the amsart/book/proc definition of \MR.
\providecommand{\MRhref}[2]{%
  \href{http://www.ams.org/mathscinet-getitem?mr=#1}{#2}
}
\providecommand{\href}[2]{#2}
\begin{thebibliography}{Wyl08}

\bibitem[Dem]{agbook}
Jean-Pierre Demailly, \emph{Complex analytic and differential geometry}.

\bibitem[SZ24]{SunZhang}
Song Sun and Junsheng Zhang, \emph{K\"ahler-{R}icci shrinkers and {F}ano fibrations}, arXiv preprint arXiv:2410.09661 (2024).

\bibitem[Wu03]{Wu}
Siye Wu, \emph{On the instanton complex of holomorphic {M}orse theory}, Comm. Anal. Geom. \textbf{11} (2003), no.~4, 775--807. \MR{2015176}

\bibitem[Wyl08]{Wylie}
William Wylie, \emph{Complete shrinking {R}icci solitons have finite fundamental group}, Proc. Amer. Math. Soc. \textbf{136} (2008), no.~5, 1803--1806. \MR{2373611}

\end{thebibliography}

\typeout{get arXiv to do 4 passes: Label(s) may have changed. Rerun}
\end{document}